\documentclass[runningheads,a4paper]{llncs}

\usepackage{amssymb}
\usepackage{graphicx,longtable,wrapfig}
\usepackage{amsmath}
\usepackage{epstopdf}

\usepackage{url}

\begin{document}
\newcommand{\E}{\mathsf{E}}
\newcommand{\Prob}{\mathsf{P}}
\newcommand{\G}{\mathsf{\Gamma}}
\newcommand{\BA}{Barab\'asi--Albert}
\newcommand{\BO}{Buckley--Osthus}
\newcommand{\HK}{Holme and Kim}

\newcommand{\widesim}[2][1.5]{
  \mathrel{\overset{#2}{\scalebox{#1}[1]{$\sim$}}}
}

\mainmatter

\title{Generalized preferential attachment: \\ tunable power-law degree distribution and clustering coefficient}

\titlerunning{Generalized preferential attachment}

\author{Liudmila~Ostroumova\inst{1}\inst{2}
\and Alexander~Ryabchenko\inst{1}\inst{3}
\and Egor~Samosvat\inst{1}\inst{3}}

\institute{Yandex, Moscow, Russia
\and
Moscow State University, Moscow, Russia
\and
Moscow Institute of Physics and Technology, Moscow, Russia}

%\author{Liudmila Ostroumova \and Alexander Ryabchenko \and Egor Samosvat}
%\authorrunning{Ostroumova L., Ryabchenko A., Samosvat E.}
%\urldef{\mailsa}\path|{ostroumova-la, d-sun-d, sameg}@yandex-team.ru|
%\institute{Yandex, Moscow, Russia\\
%\mailsa\\
%}

\maketitle

\begin{abstract}
We propose a common framework for analysis of a wide class of preferential attachment models, which includes LCD, Buckley--Osthus, Holme--Kim and many others. The class is defined in terms of constraints that are sufficient for the study of the degree distribution and the clustering coefficient. We also consider a particular parameterized model from the class and illustrate the power of our approach as follows. Applying our general results to this model, we show that both the parameter of the power-law degree distribution and the clustering coefficient can be controlled via variation of the model parameters. In particular, the model turns out to be able to reflect realistically these two quantitative characteristics of a real network, thus performing better than previous preferential attachment models. All our theoretical results are illustrated empirically.
\vspace{20pt}

\end{abstract}

\section{Introduction}\label{Introduction}

\footnotetext[4]{The authors are given in alphabetical order}

Numerous random graph models have been proposed to reflect and predict important quantitative and topological aspects of growing real-world networks, from Internet and society \cite{BA_Review,Networks,Chayes} to biological networks \cite{BioInfoPrior}. Such models are of use in experimental physics, bioinformatics, information retrieval, and data mining. An extensive review can be found elsewhere (e.g., see~\cite{BA_Review,Networks,Math_Results}). Though largely successful in capturing key qualitative properties of real-world networks, such models may lack some of their important characteristics.

The simplest characteristic of a vertex in a network is the degree, the number of adjacent edges. Probably the most extensively studied property of networks is their vertex degree distribution. For the majority of studied real-world networks, the portion of vertices with degree $d$ was observed to decrease as $d^{-\gamma}$, usually with $2 < \gamma < 3$, see \cite{BA,Networks,Broder,F-F-F}. Such networks are often called scale-free.

Another important characteristic of networks is their clustering coefficient, a measure capturing the tendency of a network to form clusters, densely interconnected sets of vertices. Various definitions of the clustering coefficient can be found in the literature, see \cite{Math_Results} for a discussion on their relationship. We consider the most popular two: the global clustering coefficient and the average local clustering coefficient (see Section~\ref{ClusteringCoefficient} for definitions). For the majority of studied real-world networks, the average local clustering coefficient varies in the range from $0.01$ to $0.8$ and does not change much as the network grows \cite{Networks}. Modeling real-world networks with accurately capturing not only their power-law degree distribution, but also clustering coefficient, has been a challenge.

In order to combine tunable degree distribution and clustering in one model, some authors \cite{BioInfoPrior,Prior1,Prior2} proposed to start with a concrete prior distribution of vertex degrees and clustering and then generate a random graph under such constraints. However, adjusting a model to a particular graph seems to be not generic enough and can be suspected in ``overfitting''. A more natural approach is to consider a graph as the result of a random process defined by certain reasonable realistic rules guaranteeing the desired properties observed in real networks. Perhaps the most widely studied realization of this approach is preferential attachment. In Section \ref{BackgroundPAmodels}, we give a background on previous studies in this field.

In this paper, we propose a new class of preferential attachment random graph models thus generalizing some previous approaches. We study this class theoretically: we prove the power law for the degree distribution and approximate the clustering coefficient. We demonstrate that in preferential attachment graphs two definitions of the clustering coefficient give quite different values. We also propose a concrete parameterized model from our class where both the power-law exponent and the clustering coefficient can be tuned. All our theoretical results are illustrated
 experimentally.

The remainder of the paper is organized as follows. In Section~\ref{BackgroundPAmodels}, we give a background on previous studies of preferential attachment models. In Section~\ref{Theory}, we propose a definition of a new class of models, and obtain some general results for all models in this class. Then, in Section~\ref{Polynom}, we describe one particular model from the proposed class. We demonstrate results obtained for graphs generated in this model in Section \ref{Experiments}. Section \ref{Conclusion} concludes the paper.

\section{Preferential Attachment Random Graph Models}
\label{BackgroundPAmodels}

In 1999, Barab\'asi and Albert observed~\cite{BA} that the degree distribution of the World Wide Web follows the power law with the parameter $\sim 2.1$. As a possible explanation for this phenomenon, they proposed a graph construction stochastic process, which is a Markov chain of graphs, governed by the \emph{preferential attachment}. At each time step in the process, a new vertex is added to the graph and is joined to $m$ different vertices already existing in the graph chosen with probabilities proportional to their degrees.

Denote by $d_v^n$ the degree of a vertex $v$ in the growing graph at time $n$. At each step $m$ edges are added, so we have $\sum_v d_v^n = 2 m n$. This observation and the preferential attachment rule imply that
\begin{equation}\label{SimplePrefAttach}
  \textbf{P}(d_v^{n+1} = d + 1 \mid  d_v^n = d) = \frac{d}{2n} \;,
\end{equation}
where $\mathbf{P}$ denotes the probability of an event. Note that the condition (\ref{SimplePrefAttach}) on the attachment probability does not specify the distribution of $m$ vertices to be joined to, in particular their dependence. Therefore, it would be more accurate to say that Barab\'asi and Albert proposed not a single model, but a class of models. As it was shown later, there is a whole range of models that fit the \BA\ description, but possess very different behavior.

\begin{theorem}[Bollob\'as, Riordan \cite{Math_Results}]
Let $f(n)$, $n \geq 2$, be any integer valued function with $f(2) = 0$ and $f(n) \leq f(n + 1) \leq f(n) + 1$ for every $n \geq 2$, such that $f(n) \rightarrow \infty$ as $n \rightarrow \infty$. Then there is a random graph process $T(n)$ satisfying (\ref{SimplePrefAttach}) such that, with probability $1$, $T(n)$ has exactly $f(n)$ triangles for all sufficiently large $n$.
\end{theorem}

In \cite{LCD_degrees}, Bollob\'as and Riordan proposed a concrete precisely defined model of the \BA\ type, known as the LCD-model, and proved that for $d < n^{\frac{1}{15}}$, the portion of vertices with degree $d$ asymptotically almost surely obeys the power law with the parameter $3$. Recently Grechnikov substantially improved this result \cite{Gr} and removed the restriction on $d$. It was shown also that the expectation of the global clustering coefficient in the model is asymptotically proportional to $\frac{(\log n)^2}{n}$ and therefore tends to zero as the graph grows \cite{Math_Results}.

One obtains a natural generalization of the LCD-model, requiring the probability of attachment of a new vertex $n+1$ to a vertex $v$ to be proportional to $d_v^n + m \beta$, where $\beta$ is a constant representing the \emph{initial attractiveness} of a vertex. Buckley and Osthus~\cite{Buckley_Osthus} proposed a precisely defined model with a nonnegative integer $\beta$. M\'ori~\cite{Mori} generalized this model to real~$\beta > -1$. For both models, the degree distribution was shown to follow the power law with the parameter $3 + \beta$ in the range of small degrees. The recent result of Eggemann and Noble~\cite{Mori_Clustering} implies that the expectation of the global clustering coefficient in the M\'ori model with $\beta > 0$ is asymptotically proportional to $\frac{\log n}{n}$. For $\beta = 0$, the M\'ori model is almost identical to the LCD-model. Therefore the authors of \cite{Mori_Clustering} emphasize the confusing difference between clustering coefficients ($ \frac{(\log n)^2}{n}$~versus~$\frac{\log n}{n}$).

The main drawback of the described preferential attachment models is unrealistic behavior of the clustering coefficient. In fact, for all discussed models the clustering coefficient tends to zero as a graph grows, while in the real-world networks the clustering coefficient is approximately a constant \cite{Networks}.

A model with asymptotically constant (average local) clustering coefficient was proposed by Holme and Kim~\cite{Holme_Kim}. The idea is to mix preferential attachment steps with the steps of triangle formation. This model allows to tune the clustering coefficient by varying the probability of the triangle formation step. However, experiments and empirical analysis show that the degree distribution in this model obeys the power law with the fixed parameter close to $3$, which does not suit most real networks. RAN (random Apollonian network) proposed in \cite{RAN} is another interesting example of a Barab\'asi-Albert type model with asymptotically constant (average local) clustering.

There is a variety of other models, not mentioned here, that are also based on the idea of preferential attachment. Analyses of properties of all these models are often very similar. In the next section, we consider theorems aimed at simplifying these analyses and providing a general framework for them. In order to do this, we define a new class of preferential attachment models that generalizes models mentioned above, as well as many others. We also propose a new parameterized model which belongs to this class that allows to tune both the power-law exponent and the clustering coefficient by adjusting the parameters.

\section{Theoretical Results}\label{Theory}

In this section, we define a general class of preferential attachment models. For all models in this class we are able to prove the power-law degree distribution.
If an additional property is fulfilled, we are able to analyze the behavior of the clustering coefficient as the network grows.
%We also estimate the number of pairs of adjacent edges in models from this class and therefore can analyze the behavior of the clustering coefficient as the network grows.

\subsection{Definition of the $PA$-class}\label{class}

Let $G_{m}^n$ ($n \ge n_0$) be a graph with $n$ vertices $\{1, \ldots, n\}$ and $mn$ edges obtained as a result of the following random graph process. We start at the time $n_0$ from an arbitrary graph $G_{m}^{n_0}$ with $n_0$ vertices and $m n_0$ edges. On the $(n+1)$-th step ($n\geq n_0$), we make the graph $G_{m}^{n+1}$ from $G_{m}^{n}$ by adding a new vertex $n+1$ and $m$ edges connecting this vertex to some $m$ vertices from the set $\{1, \ldots , n, n+1\}$. Denote by $d_v^{n}$ the degree of a vertex $v$ in~$G_m^n$. If for some constants $A$ and $B$ the following conditions are satisfied
\begin{equation}\label{OneStepChangedDegreeDraft}
\textbf{P}\left( d_v^{n+1} = d_v^{n} \mid G_m^{n}\right) = 1 - A \frac{d_v^{n}}{n} - B\frac{1}{n} + O\left(\frac{\left(d_v^n\right)^2}{n^2}\right),\,\,1 \le v \le n \;,
\end{equation}
\begin{equation}\label{OneStepChangedDegreeDraft_2}
\textbf{P}\left( d_v^{n+1} = d_v^{n} + 1 \mid G_m^{n}\right) =  A \frac{d_v^{n}}{n} + B\frac{1}{n} + O\left(\frac{\left(d_v^n\right)^2}{n^2}\right), \,\,1 \le v \le n \;,
\end{equation}
\begin{equation}\label{OneStepChangedDegreeDraft_3}
\textbf{P}\left( d_v^{n+1} = d_v^{n} + j \mid G_m^{n}\right) =  O\left(\frac{\left(d_v^n\right)^2}{n^2}\right), \,\,
2\le j \le m,\,\,1 \le v \le n \;,
\end{equation}
\begin{equation}\label{LoopProbability}
\textbf{P}( d_{n+1}^{n+1} =  m + j ) = O\left(\frac 1 n \right), \,\,
1\le j \le m\;,
\end{equation}
then we say that the random graph process $G_m^n$ is a model from the $PA$-class. Condition~(\ref{LoopProbability}) means that the probability to have a self-loop in the added vertex is small. As we will show later, certain minor details of the models from this class, such as whether loops and multiple edges are allowed, are irrelevant.

Since we add $m$ edges at each step, summing up the equalities \eqref{OneStepChangedDegreeDraft_2}-\eqref{LoopProbability} (with corresponding coefficients) over all vertices and neglecting error terms we get $2mA + B = m$. It is possible to prove that the sum of error terms in this case is $0$, but for simplicity we just set $2mA + B = m$. Furthermore,
we have $0 \le A \le 1$ (for \eqref{OneStepChangedDegreeDraft_2} we need $mA + B \ge 0$ and we set  $2mA + B = m$, therefore $A \le 1$).

%since the probabilities defined in (\ref{OneStepChangedDegreeDraft}) and (\ref{OneStepChangedDegreeDraft_2})
%must be positive for all $d_i^n \ge m$ and~all~$n$ (e.g., for large $n$ and $d_i^n = n$ these probabilities tends to $1-A$ and $A$ respectively).
%Note that if $A = 1$ then $B = -m$ and probability of attachment to a node is proportional to its incoming degree.

Here we want to emphasize that we indeed defined not a single model but a class of models. Even fixing values of parameters $A$ and $m$ does not specify a concrete procedure for constructing a network. What this definition lacks is the precise description of the distribution of vertices a new incoming vertex is being connected to, and therefore there is a range of models possessing very different properties and satisfying the conditions~(\ref{OneStepChangedDegreeDraft}--\ref{LoopProbability}). For example, the LCD, the Holme--Kim and the RAN models belong to the $PA$-class with $A = 1/2$ and $B = 0$. The Buckley--Osthus (M\'ori) model also belongs to the $PA$-class with $A = \frac{1}{2+\beta}$ and $B = \frac{m \beta}{2+\beta}$. Another example is considered in detail in Sections~\ref{Polynom} and~\ref{Experiments}. This situation is somewhat similar to that with the definition of the \BA\ models, though our class is wider in a sense that the exponent of the power-law degree distribution is tunable.

In mathematical analysis of network models, there is a tendency to consider only fully and precisely defined models. In contrast, we provide results about general properties for the whole $PA$-class in the next two subsections.

\subsection{Power Law Degree Distribution}\label{DegreeDistribution}

Even though the precise description of the distribution of vertices a new incoming vertex is going to be connected to is not specified, we are still able to describe the degree distribution of the network.

First, we estimate $N_n(d)$, the number of vertices with given degree $d$ in~$G_m^n$. We prove the following result on the expectation $\E N_n(d)$ of $N_n(d)$.

\begin{theorem}\label{Expectation}
For every $d \ge m$ we have
$\E N_n(d) = c(m,d) \left(n + O\left(d^{2 + \frac{1}{A}}\right)\right)$,
where
\begin{equation*}\label{Constant}
c(m,d) = \frac{\G\left(d + \frac{B}{A}\right)\G\left(m + \frac{B+1}{A}\right)}{A\G\left(d + \frac{B+A+1}{A}\right)\G\left(m + \frac{B}{A}\right)}
\widesim{d \rightarrow \infty} \frac{\G\left(m + \frac{B+1}{A}\right)d^{-1-\frac{1}{A}}}{A \G\left(m + \frac{B}{A}\right)}  \;,
\end{equation*}
and $\G(x)$ is the gamma function.
\end{theorem}

Second, we show that the number of vertices with given degree $d$ is highly concentrated around its expectation.
\begin{theorem}
\label{Concentration}
For every model from the PA-class and for every $d=d(n)$ we have
$$
\Prob\left(|N_n(d) - \E N_n(d)| \ge d \, \sqrt{n}  \, \log{n}\right) = O\left(n^{-\log{n}}\right).
$$
Therefore, for any $\delta > 0$ there exists a function $\varphi(n) \in o(1)$ such that
$$
\lim_{n \to \infty} \Prob \left(\exists \, d \le n^{\frac{A-\delta}{4A+2}}:|N_n(d) - \E N_n(d)| \ge \varphi(n)\,\E N_n(d)\right) = 0 \;.
$$
\end{theorem}
These two theorems mean that the degree distribution follows
(asymptotically) the power law with the parameter $1+\frac{1}{A}$.

Theorem \ref{Expectation} is proved by induction on $d$ and $n$. It is easy to see that given a graph $G_m^n$, we can express the conditional expectation of the number of vertices with degree $d$ in $G_m^{n+1}$ (i.e., $\E ( N_{n+1}(d)\mid G_m^n)$) in terms of $N_n(d), N_n(d-1), \dots, N_n(d-m)$. Here we use the fact that the probability of having an edge between the vertex $n+1$ and a vertex $v$ depends on the degree of $v$ (see (\ref{OneStepChangedDegreeDraft})). Using the law of total expectation we obtain the recurrent relation for $ \E N_{n+1}(d)$ and prove the statement of Theorem \ref{Expectation} by induction.

We use the Azuma--Hoeffding inequality to prove the concentration result of Theorem~\ref{Concentration}. In order to do this, we consider the martingale $X_i(d) = \E(N_n(d)\mid G_m^i)$, $i = 0, \ldots, n$. The complete proofs of these theorems are technical and are placed in Appendix due to space constraints.

\subsection{Clustering Coefficient}\label{ClusteringCoefficient}

Here we consider the clustering coefficient in models of the $PA$-class.  There are two popular definitions of the clustering coefficient. The \emph{global clustering coefficient} $C_1(n)$ is the ratio of three times the number of triangles to the number of pairs of adjacent edges in G. The \emph{average local clustering coefficient} is defined as follows: $C_2(n) = \frac{1}{n} \sum_{i=1}^n C(i)$, where $C(i)$ is the local clustering coefficient for a vertex $i$: $C(i) = \frac{T^i}{P_2^i}$, where $T^i$ is the number of edges between neighbors of the vertex $i$ and $P_2^i$ is the number of pairs of neighbors. Results for some classical preferential attachment models (LCD and M\'ori) are mentioned in Section~\ref{BackgroundPAmodels}.

Here we generalize these results. First, we study the random variable $P_2(n)$ equal to the number of $P_2$'s in a random graph $G_m^n$ from an arbitrary model that belongs to the PA-class. In the theorems below, we use the following notation. By {\bf whp} (``with high probability'') we mean that for some sequence ${A_n}$ of events, $P(A_n) \to 1$ as $n \to \infty$. We say $a_n \sim b_n$ if $a_n = (1+o(1))b_n$, and we say $a_n \propto b_n$ if $C_0 b_n \le a_n \le C_1 b_n$ for some constants $C_0, C_1 > 0$.

\begin{theorem}\label{P_2}
For every model from the $PA$-class, we have
\begin{itemize}
\item[(1)] if $2A<1$, then \textbf{whp}
$P_2(n) \sim \left(2m(A+B) + \frac{m(m-1)}{2}\right) \frac{n}{1-2A} \;,$
\item[(2)] if $2A=1$, then \textbf{whp} $P_2(n) \sim \left(2m(A+B) + \frac{m(m-1)}{2}\right)
n \log(n) \;,$
\item[(3)] if $2A>1$, then for any $\varepsilon>0$ \textbf{whp} $n^{2A-\varepsilon} \le P_2(n) \le n^{2A+\varepsilon}$.
\end{itemize}
\end{theorem}

The ideas of the proof of Theorem~\ref{P_2} are given in Appendix. Here it is worth noting that the value $P_2(n)$ in scale-free graphs is usually determined by the power-law exponent $\gamma$. Indeed, we have
$
P_2(n) = \sum_{d=1}^{d_{\max}} N_n(d)\frac{d(d-1)}{2}
\propto \sum_{d=1}^{d_{\max}} n d^{2-\gamma},
$
where  $d_{\max}$ is the maximum possible degree of a vertex in $G_m^n$. Therefore if $\gamma > 3$, then $P_2(n)$ is linear in $n$. However, if $\gamma \le 3$, then $P_2(n)$ is superlinear.

Next, we study the random variable $T(n)$ equal to the number of triangles in $G_m^n$. Note that in any model from the $PA$-class we have $T(n) = O(n)$ since at each step we add at most $\frac{m(m-1)}{2}$ triangles. If we combine this fact with the previous observation, we see that if $\gamma \leq 3$, then in any preferential attachment model (in which the out-degree of each vertex equals $m$) the global clustering coefficient tends to zero as $n$ grows.

Our aim is to find models with constant clustering coefficient. Let us consider a subclass of the $PA$-class with the following property:
\begin{equation}\label{D_definition}
\textbf{P}\left( d_i^{n+1} = d_i^{n} + 1, d_j^{n+1} = d_j^{n} + 1 \mid G_m^{n}\right) = e_{ij} \frac{D}{mn} + O\left(\frac{d_i^{n} d_j^{n}}{n^2}\right) \;.
\end{equation}
Here $e_{ij}$ is the number of edges between vertices $i$ and $j$ in $G_m^n$ and $D$ is a positive constant. Note that this property still does not define the correlation between edges completely.

\begin{theorem}\label{Triangles}
Let $G_m^n$ satisfy the condition (\ref{D_definition}). Then {\bf whp}
$
T(n) \sim D \, n \;.
$
\end{theorem}
The proof of this theorem is straightforward. The expectation of the number of triangles we add at each step is $D + o(1)$. The fact that the sum of $O\left(\frac{d_i^{n} d_j^{n}}{n^2}\right)$ over all \textit{adjacent} vertices is $o(1)$ can be shown by induction using the conditions (\ref{OneStepChangedDegreeDraft}--\ref{LoopProbability}). It is also possible to first prove that the maximum degree grows as $n^A$ and then use this fact to estimate the sum of error terms. Therefore $\E T(n) = D n + o(n)$. The Azuma--Hoeffding inequality can be used to prove concentration.

As a consequence of Theorems \ref{P_2} and \ref{Triangles}, we get the following result on the global clustering coefficient $C_1(n)$ of the graph $G_m^n$.
\begin{theorem}\label{Cluster}
Let $G_m^n$ belong to the $PA$-class and satisfy the condition (\ref{D_definition}). Then
\begin{itemize}
\item[(1)]
If $2A<1$ then \textbf{whp}
$C_1(n) \sim  \frac{3(1-2A)D}{\left(2m(A+B) + \frac{m(m-1)}{2}\right)} \;,$
\item[(2)]
If $2A=1$ then \textbf{whp}
$C_1(n) %\sim \frac{1}{\log{n}}\;.$
\sim  \frac{3 D}{\left(2m(A+B) + \frac{m(m-1)}{2}\right) \log n} \;,$ % = \frac{6D}{\log(n) m (m + 1)}
\item[(3)]
If $2A>1$ then for any $\varepsilon > 0$ \textbf{whp}
$n^{1-2A-\varepsilon} \le C_1(n) \le n^{1-2A+\varepsilon}\;.$
%$C_1(n) \propto n^{1-2A}\;.$
\end{itemize}
\end{theorem}

Theorem \ref{Cluster} shows that in some cases ($2A\ge1$) the global clustering coefficient $C_1(n)$ tends to zero as the number of vertices grows. We empirically show in Section \ref{Experiments} that the average local clustering coefficient $C_2(n)$ behaves differently.

The theoretical analysis in this case is much harder, but we can easily show why $C_2(n)$ does not tend to zero if the condition (\ref{D_definition}) holds. From Theorems \ref{Expectation} and \ref{Concentration} it follows that {\bf whp} the number of vertices with degree $m$ in $G_m^n$ is greater than $c n$ for some positive constant $c$. The expectation of the number of triangles we add at each step is $D + o(1)$. Therefore {\bf whp}
$
C_2(n) \ge \frac{1}{n} \sum_{i: \deg(i) = m} C(i) \ge \frac{2 c D}{m(m+1)} \, .
$

In the next section we introduce a concrete nontrivial model from the $PA$-class.

\section{Polynomial Model}\label{Polynom}

In this section, we consider {\it polynomial random graph models} that belong to the general $PA$-class defined above.
 %In Subsection~\ref{PolynomDef} we give definitions. In Subsection~\ref{PropertiesPolynom}, we find the relations between the polynomial model parameters and those of the $PA$-class.
Applying our theoretical results to polynomial models, we find the model to be very flexible: one can tune the parameter of the degree distribution and the clustering coefficient.

\vspace{10px}

\textbf{Definition of Polynomial Model } %\label{PolynomDef}
Let us define the \textit{polynomial model}. As in the random graph process from Subsection~\ref{class}, we construct a graph $G_m^n$ step by step. On the $(n+1)$-th step the graph $G_m^{n+1}$ is made from the graph $G_m^n$ by adding a new vertex $n + 1$ and sequentially drawing $m$ edges (multiple edges and self-loops are allowed).
%$e_1, \ldots, e_m$ connecting this vertex to some $m$ vertices $i_1, \dots, i_m \in \{1, \dots, n+1\}$. Some of $i_1, \dots, i_m$ can be equal to each other, %so multiple edges are allowed. Also $i_j = n+1$ for some $j$ results in a self-loop (multiple self-loops are also allowed).

We say that an edge $ij$ is directed from $i$ to $j$ if $i \ge j$, so the out-degree of each vertex equals $m$. We also say that $i$ and $j$ are respectively \textsl{source} and \textsl{target} ends of $ij$. We consider different approaches to add new edges from the vertex $n + 1$. We first choose an edge from the existing graph $G_m^n$ uniformly and independently at random and then have three options:
\begin{itemize}
  \item Preferential attachment (PA): draw one edge from $n+1$ to the \emph{target} end of the chosen edge
  \item Uniform (U): draw one edge from $n+1$ to the \emph{source} end of the chosen edge
  \item Triangle formation (TF): draw two edges from $n+1$ to \emph{target} and \emph{source} ends of the chosen edge
\end{itemize}

Let us now specify how to draw $m$ edges from the vertex $n+1$. Consider a collection of positive parameters $\{\alpha_{k,l}\}$ for ${0 \le k \le m/2}$
and ${0 \le l \le m - 2k}$ such that ${\sum_{k,l} \alpha_{k,l} = 1}$, these parameters fully define our model. At the beginning of the $n+1$ step with probabilities $\{\alpha_{k,l}\}$ we choose some $k = k_0$ and $l = l_0$, then we draw $l_0$ edges using PA, $2 k_0$ edges using TF and ${(m - l_0 - 2k_0)}$ edges using~U.
This random graph process defines the polynomial model and from the definition it follows that graphs in this model
can be generated in linear time. This model belongs to the $PA$-class.
Indeed, one can formally show by simple calculations that the conditions (\ref{OneStepChangedDegreeDraft}--\ref{LoopProbability}) hold for this model.

At this point the model is defined but let us explain why we call it polynomial.
Denote by $\widehat{d^{n}_{i}} = d^{n}_{i} - m$  the in-degree  of a vertex $i$ in $G_m^n$. Let us recall that by $e_{ij}$ we denote the number of edges between vertices $i$ and $j$. For every $k,l$ such that $0 \le k \le m/2$ and $0 \le l \le m - 2k$, let
$
M_{k,l}^{n,m}(i_1,\dots,i_m) = \frac{1}{n^{m-l-2k}} \prod_{x=1}^k \frac{e_{i_{2x} i_{2x-1}}}{2 m n} \prod_{y=2k+1}^{2k+l} \frac{\widehat{d^{n}_{i_y}}}{m n}.
$
This is a monomial depending on $\widehat{d^{n}_{i_y}}$ and $e_{i_{2x} i_{2x-1}}$.  We define the polynomial $\sum_{k,l} \alpha_{k,l} M_{k,l}^{n,m}(i_1,\dots,i_m)$. It is easy to check that
\begin{multline}\label{PolinomialModel}
\mathbf{P} \left(\text{edges $e_1, \ldots, e_m$ go to vertices $i_1, \ldots , i_m$, respectively} \right) = \\
=\sum_{k=0}^{m/2}  \sum_{l=0}^{m-2k}
 \alpha_{k,l} M_{k,l}^{n,m}(i_1,\dots,i_m) \;.
\end{multline}
%Also $(d^{n}_{i_y})^{in}$ can be replaced by $d^{n}_{i_y}+m$  and polynomial can be rewritten in terms~of~$d_i^n$.

Many models are special cases of the polynomial model.
If we consider the polynomial $\prod_{y=1}^{m} \frac{\widehat{d^{n}_{i_y}} + m}{2mn}$, then we obtain a model that is practically identical to the LCD-model. The Buckley--Osthus model can be also interpreted in terms of the polynomial model.

\vspace{10px}

\textbf{Properties }%\label{PropertiesPolynom}
It is easy to check that the parameters $\alpha_{k,l}$ from (\ref{PolinomialModel}) and $A$ from (\ref{OneStepChangedDegreeDraft}) are related in the following way:
\begin{gather}\label{Parameters}
A = \sum \alpha_{k,l} \frac{l+k}{m} \;. %\\
\end{gather}
This means that we can use an arbitrary value of $A \in [0,1]$ and any power-law exponent $\gamma \in (2, \infty)$ in the graph generation.
Also note that $D = \sum_{k,l} k \alpha_{k,l}$\ .

In the next section we analyze experimentally some properties of graphs in the polynomial model.
We generate polynomial graphs and compare their properties with theoretical results we obtained.

\section{Experiments}\label{Experiments}

In this section, we choose a three-parameter model from the family of polynomial graph models defined in Section \ref{Polynom} and analyze the properties of the generated graphs depending on the parameters.

\subsection{Description of Empirically Studied Polynomial Model}

We study empirically graphs in the polynomial model with $m = 2p$
and the probability to draw edges to vertices $i_1, \ldots, i_{2m}$ equals
\begin{equation*}
\prod_{k=1}^{p}\left(\alpha \frac{\widehat{d^{n}_{i_{2k}}}\widehat{d^{n}_{i_{2k-1}}}}{(mn)^2} + \beta\frac{ e_{i_{2k} i_{2k-1}}}{2mn} + \frac{\delta}{(n)^2} \right).
\end{equation*}
Here we need $\alpha, \beta, \delta \ge 0$  and $\alpha + \beta + \delta = 1$, therefore, we have three independent model parameters: $m$, $\alpha$, and~$\beta$. Note that here we write the polynomial in a symmetric form as we ignore the order of edges.

Based on our theoretical results, we have certain expectations about the properties of generated graphs. From (\ref{Parameters}) we obtain that in this model $A = \alpha + \frac{\beta}{2}, B = m(\delta - \alpha)$, $D = p \beta = \frac{ m \beta}{2}$, therefore, due to Theorem \ref{Expectation} and Theorem \ref{Cluster}, we get that
\begin{equation}\label{ClustAlpha}
C_1(n) \sim  \frac{3(1 - 2\alpha - \beta)\beta}{5m - 1 - 2(2m-1)(2\alpha + \beta)}  \text{, }  \gamma = 1 + \frac{2}{2\alpha + \beta}\;.
\end{equation}

\subsection{Empirical Results}

\begin{figure}
\begin{center}
\includegraphics[height = 3.8cm]{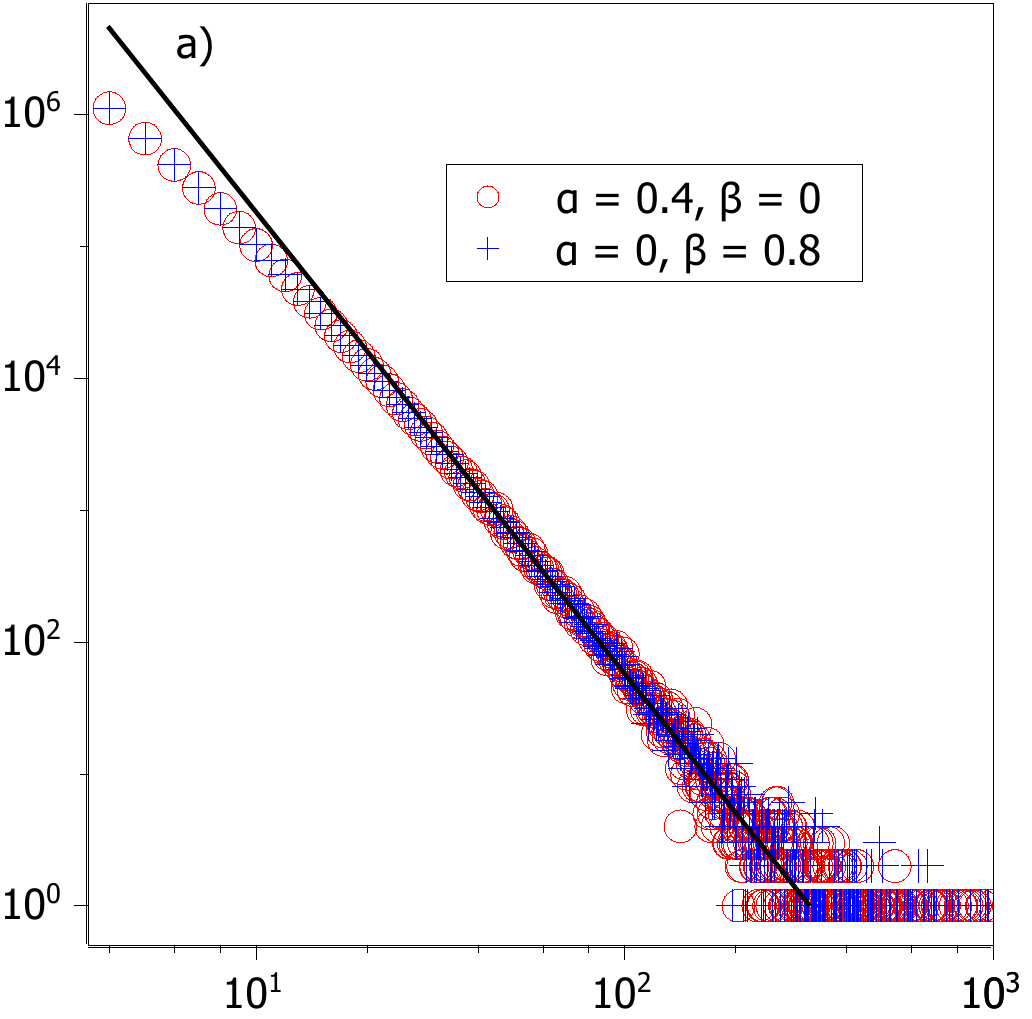}
\includegraphics[height = 3.8cm]{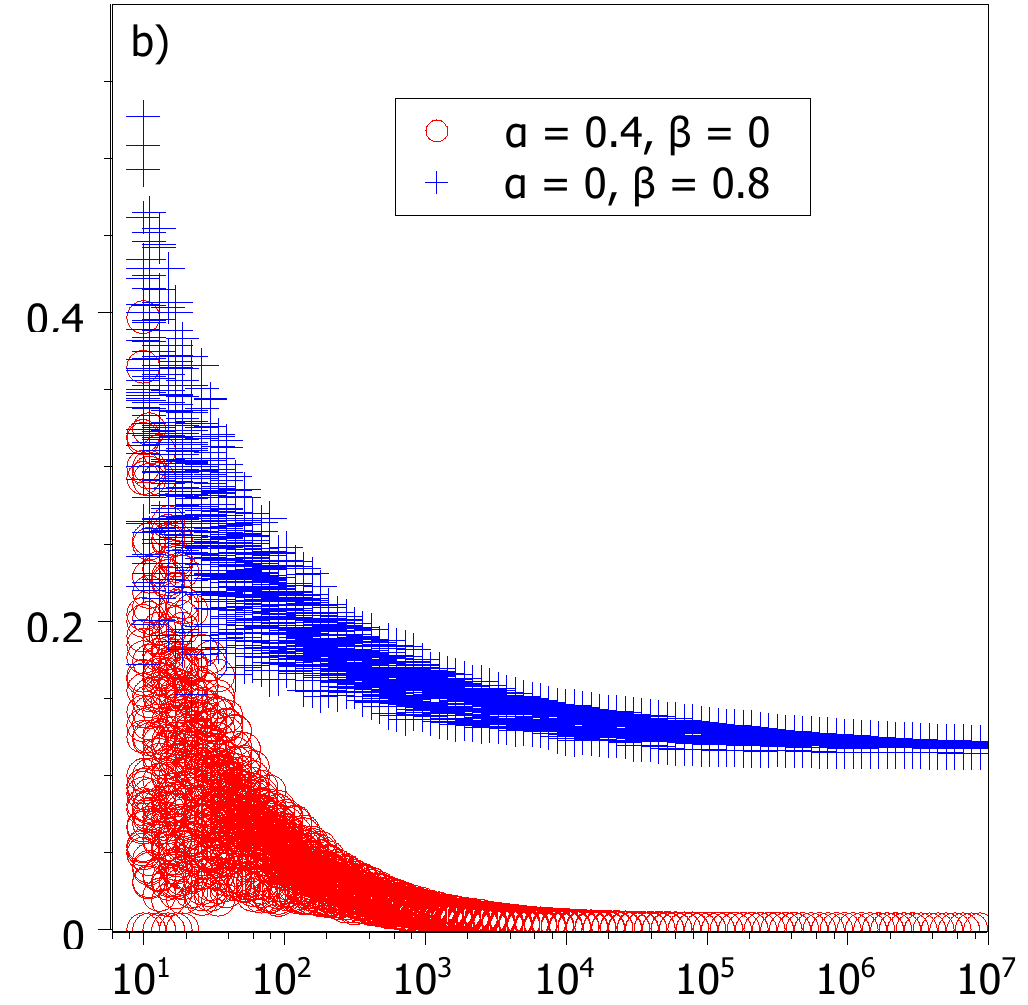}
\includegraphics[height = 3.8cm]{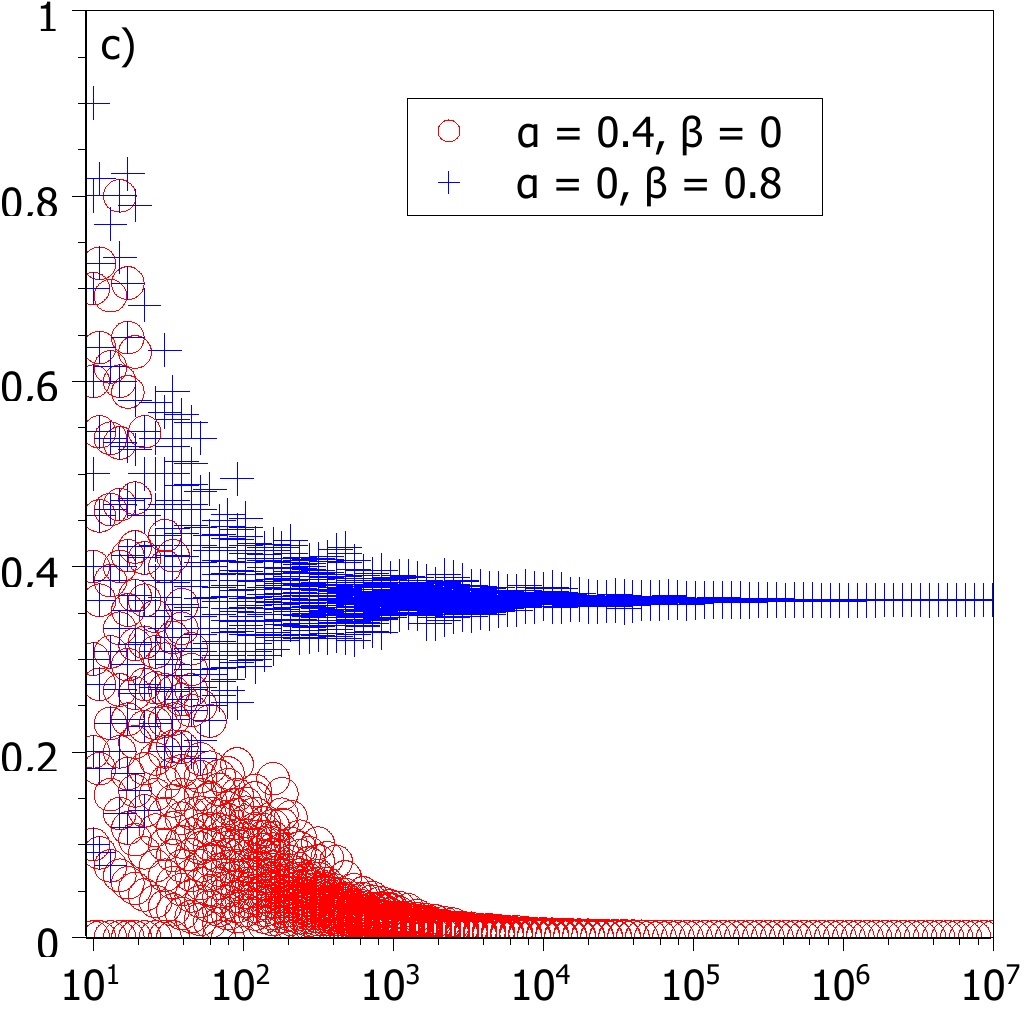}
\end{center}
\caption{ a) The degree distribution of polynomial graphs with $n = 10^7$ and $m = 2$.
b)~The global clustering coefficient of polynomial graphs with $m = 2$ depending on $n$.
c)~The average local clustering coefficient of polynomial graphs with $m = 2$ depending on $n$.}
%\textbf{[third picture has b) on it instead of c)]}
\label{fig:DegreeClust}
\end{figure}

\textbf{Degree Distribution and Clustering Coefficient }
First, we study two polynomial graphs with $n = 10^7$, $m = 2$, and $A = 0.4$, assigning $\alpha = 0.4, \beta = 0$ for the first graph and $\alpha = 0, \beta = 0.8$ for the second one. The observed degree distributions are almost identical and follow the power law with the expected parameter $\gamma = 3.5$, see Fig.~\ref{fig:DegreeClust}a.

For both cases, we also study the behavior of the global and the average local clustering coefficients of generated graphs, $40$ samples for each $n = \left[10^{1 + 0.06i}\right], i = 0,\dots,100$, see Fig.~\ref{fig:DegreeClust}bc. In the first case we observe $C_1(n) \to 0$, $C_2(n) \to 0$ (as $\beta = 0$) and in the second case $C_1(n) \to \frac{2}{15}$ (as was expected due to (\ref{ClustAlpha})) and
$C_2(n) \to \mathrm{const} > 0$.

\begin{figure}
\begin{center}
\includegraphics[height = 3.9cm]{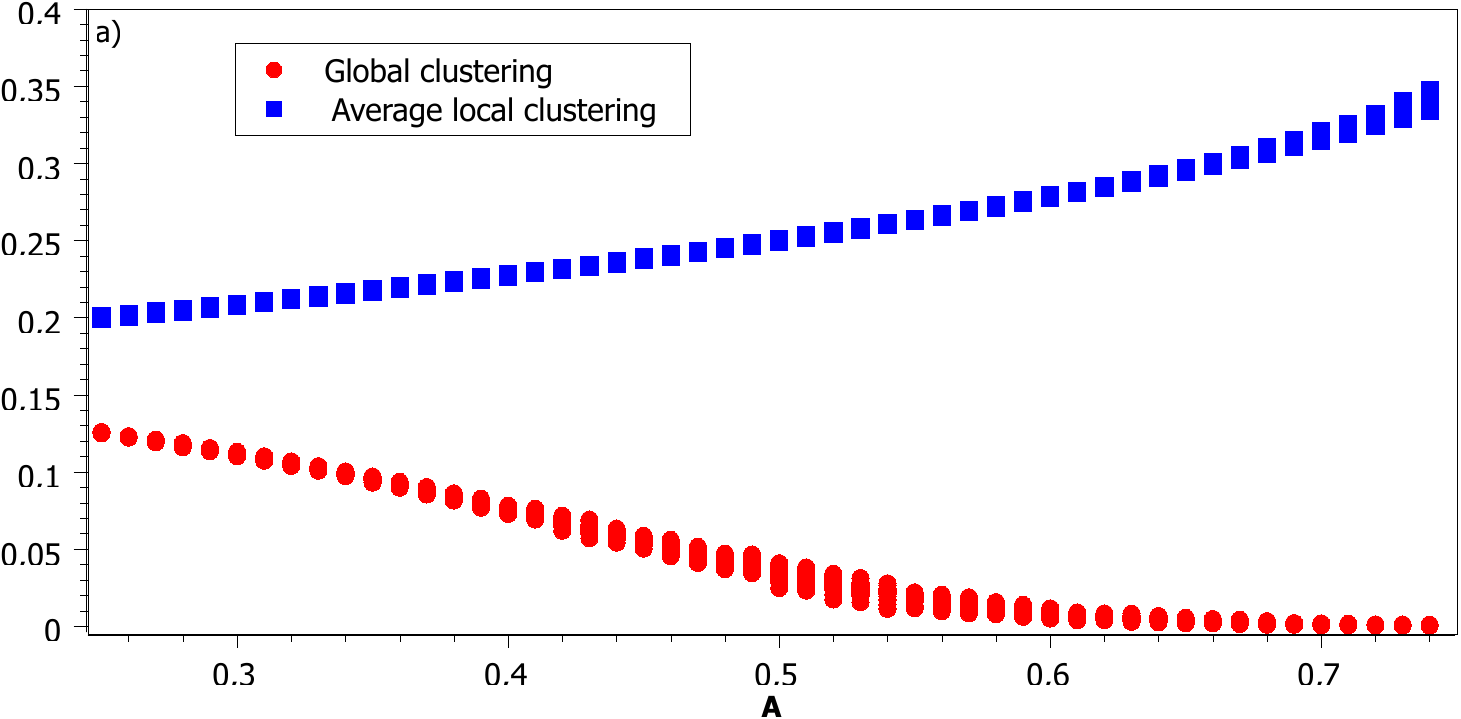}
\includegraphics[height = 3.9cm]{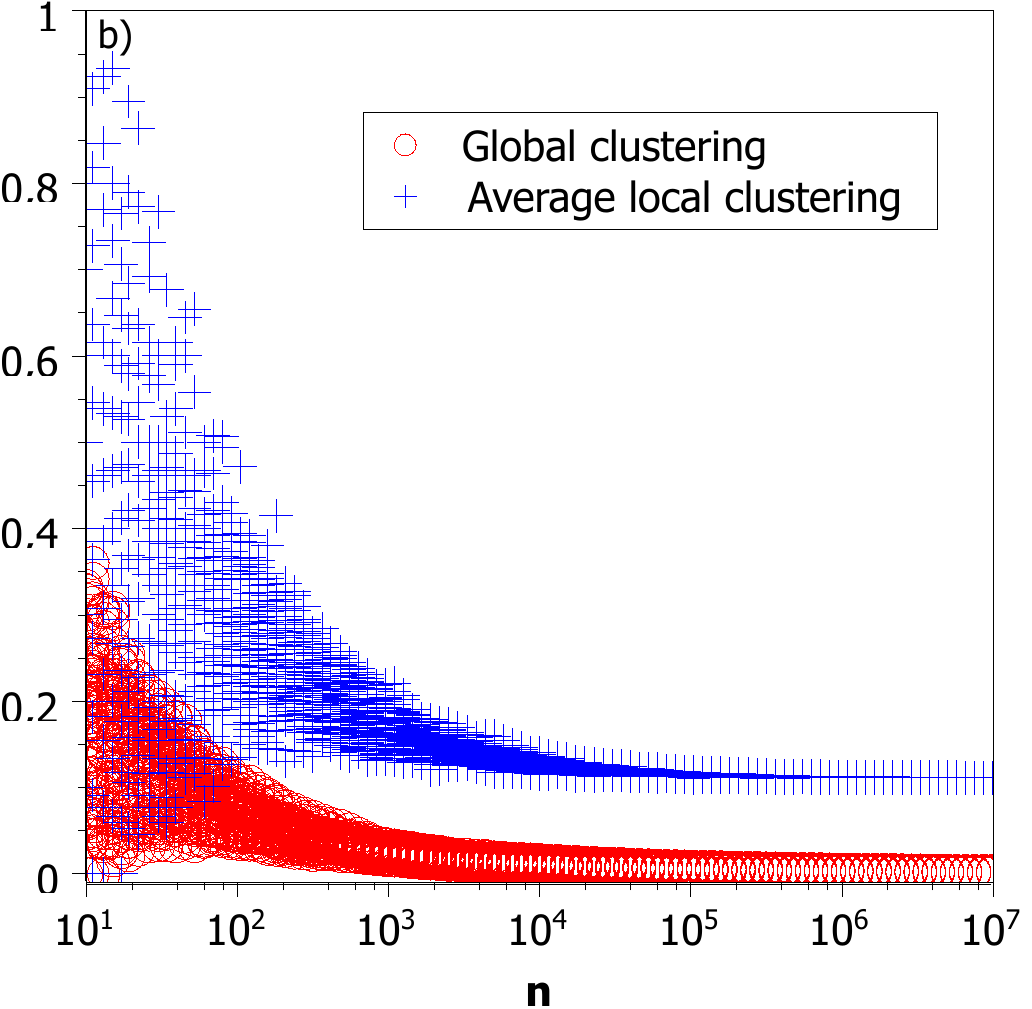}
\end{center}
\caption{a) Average local and global clustering in polynomial graphs with $n = 10^6$, $m = 2$, $\beta = 0.5$ depending on A.
b) The global and the average local clustering coefficients of polynomial graphs with $m = 2$, $\alpha = 0.5$, $\beta = 0.2$ depending on $n$.}
\label{fig:SecondFigure}
\end{figure}

We also generate graphs with $n = 10^6$, $m = 2$, and varying $A$ (we took $\beta = 0.5$ and $\alpha \in (0,0.5)$). In other words, we fix the probability of a triangle formation and vary the parameter of the power-law degree distribution. The obtained results are shown in Fig.~\ref{fig:SecondFigure}a. The behavior of the clustering coefficients is quite different. If $A$ grows, then $P_2(n)$ grows (therefore $C_1(n) \to 0$), the number of vertices with small degrees and hence high local clustering also grows (therefore $C_2(n)$ increases).

To demonstrate the difference between the global clustering and the average local clustering we generated graphs with $m = 2$, $\alpha = 0.5$, $\beta = 0.2$ and varying $n$ (Fig.~\ref{fig:SecondFigure}b). In this case we have $A = \alpha + \frac{\beta}{2} > 0.5$ and $C_1(n) \to 0$, as expected. However, for the local clustering we obtain $C_2(n) \to \mathrm{const} > 0$.

\vspace{10px}

\textbf{Comparison with Other Models }
The following table summarizes our results for the polynomial model in comparison with other mentioned preferential attachment models:

\begin{small}
\begin{center}
\begin{tabular}{|l|c|c|c|c|c|}
\hline
   &$A$ & $D$ & $\gamma$ & Global clustering & Average local clustering \\
\hline
LCD
& $1 /2 $
& $0$
& $3$
& tends to zero
& tends to zero \\
\hline
BO/M\'ori
& $1/(2+\beta)$
& $0$
& $(2,\infty)$
& tends to zero
& tends to zero\\
\hline
HK
& $1/2$
& $m_t$
& $3$
& tends to zero
& constant \\
\hline
RAN
& $1/2$
& $3$
& $3$
& tends to zero
& constant \\
\hline
Polynomial
& $\sum \alpha_{k,l} \frac{l+k}{m}$
& $\sum k \alpha_{k,l}$
& $(2,\infty)$
&constant for $A < \frac{1}{2}$
&constant \\
\hline
\end{tabular}
\end{center}
\end{small}

The polynomial model seems to be the only model where one can control the exponent in the power law of the degree distribution, and at the same time guarantee a positive clustering coefficient.

\section{Conclusions}\label{Conclusion}

In this paper, we introduced the $PA$-class of random graph models that generalizes previous preferential attachment approaches. We proved that any model from the $PA$-class possesses the power-law degree distribution with tunable parameter. We also estimated its clustering coefficient. Next, we described one particular model from the proposed class (with tunable both the degree distribution parameter and the clustering coefficient). Experiments with generated graphs illustrated our theoretical results. We also demonstrated different behavior of two versions of the clustering coefficient in preferential attachment models.

As the degree distribution of a preferential attachment model allows adjustment to reality, the clustering coefficient still gives rise to a problem in some cases. For most real-world networks the parameter $\gamma$ of their degree distribution belongs to $[2, 3]$. As we showed in Section \ref{Theory}, once $\gamma \leq 3$ in a preferential attachment model, the global clustering coefficients decreases as the graph grows, which does not correspond to the majority of real-world networks. The reason is that the number of edges added with a new vertex at each step is a constant and consequently the number of triangles grows too slowly.

Fortunately, there are many ways to overcome this obstacle. Cooper proposed a model in  which the number of added edges is a random variable \cite{CooperA1A2}. In collaboration with Pra\l at he also considered a modification of the \BA\ model, where a new vertex added at time $t$ generates $t^c$ edges \cite{Pralat}. Preferential attachment models with random initial degrees were considered in \cite{Ark}. Also there are models with adding edges between already existing nodes (e.g. \cite{CooperFrieze}). Using one of these ideas for the $PA$-class is a topic for future research.

\vspace{10px}

\textbf{Acknowledgements} Special thanks to Evgeniy Grechnikov, Gleb Gusev, Andrei Raigorodskii and anonymous reviewers for the careful reading and useful comments.

\section*{Appendix: Proofs}

\subsection*{Proof of Theorem \ref{Expectation}}

In this proof we use the notation $\theta(\cdot)$ for error terms.
By $\theta(X)$ we denote a function such that $|\theta(X)| < X$.
We also need the following notation:
\begin{equation}\label{OneStepUnChangedDegree}
\textbf{P}\left( d_v^{n+1} = d \mid d_v^{n} =d \right) = 1 - A \frac{d}{n} - B\frac{1}{n} + O\left(\frac{d^2}{n^2}\right) \; ,
\end{equation}
\begin{equation}\label{OneStepChangedDegree}
p_n^1(d) := \textbf{P}\left( d_v^{n+1} = d+1 \mid d_v^{n} = d \right) =  A \frac{d}{n} + B\frac{1}{n} + O\left(\frac{d^2}{n^2}\right) \;,
\end{equation}
\begin{equation}\label{OneStepBigChangedDegree}
p_n^j(d) := \textbf{P}\left( d_v^{n+1} = d+j \mid d_v^{n} = d \right) =  O\left(\frac{d^2}{n^2}\right),  \,\,
2\le j \le m \;.
\end{equation}
\begin{equation}%\label{LoopProbability}
p_n := \sum_{k=1}^m \textbf{P}( d_{n+1}^{n+1} =  m + k ) = O\left(\frac 1 n \right) \;.
\end{equation}
Note that the remainder term of $p_n^j(d)$ can depend on $v$.
We omit $v$ in notation $p_n^j(d)$ for simplicity of proofs.

Put $p_v(d) = \sum_{j=1}^m p_v^j(d)$.
Note that $\frac{Ad+B+1}{Ad-A+B}p_{v}^1(d-1) - p_v(d) = \frac{1}{v} + O\left(\frac{d^2}{v^2}\right)$.
%Note that in classical preferential attachment models this equality holds.
We use this equality several times in this proof.

We want to prove that $\E N_n(d) = c(m,d) \left(n + \theta\left(Cd^{2 + \frac{1}{A}}\right)\right)$ with some constant $C$ and some function $\theta$.
The proof is by induction on $d$ and then on $i$.
First, we prove the theorem for $d=m$ and all $i$. Then, if we proved the theorem for some $d=d_0$ and all $i$, we are able to prove it for $d=d_0+1$ and for all $i$.

We use the following equalities
\begin{equation}\label{d=m}
\E ( N_{i+1}(m) \mid N_i(m) ) =
 N_{i}(m) \left(1 - p_i(m)\right) + 1 - p_i\;,
\end{equation}
\begin{multline}\label{d>m}
\E ( N_{i+1}(d) \mid N_i(d), N_i(d-1), \dots, N_i(d-m) ) =
N_i(d)\left(1 - p_i(d)\right) + \\ + N_i(d-1) p_i^1(d-1) + \sum_{j=2}^m N_i(d-j) p_i^j(d-j) + O(p_i) \;.
\end{multline}

Consider the case $d=m$. %We have $p_i(m) = \frac{mA+B}{i}  + O\left(\frac{1}{i^2}\right)$.
For constant number of small $i$ we obviously have $\E N_i(m) =\frac {i} {Am+B+1} + \theta(C_1)$ with some $C_1$.
Assume that $\E N_i(m) = \frac {i} {Am+B+1} + \theta(C_1)$.
From (\ref{d=m}) we obtain
\begin{multline*}
\E N_{i+1}(m) =
\E N_{i}(m) \left(1 - p_i(m)\right) + 1 - p_i = \\
= \left( \frac {i} {Am+B+1}  + \theta(C_1)  \right)\left(1 - p_i(m)\right) + 1 + \theta(C_2/i) = \\
= \frac{i+1}{Am+B+1} + \theta(C_1)\left(1 - p_i(m)\right)
+ \theta \left(\frac{C_3}{i}\right)\frac {1} {Am+B+1}  + \theta(C_2/i) \;.
\end{multline*}
It remains to show that
$$
C_1 p_i(m) \ge \frac{C_3}{i(Am+B+1)} + \theta(C_2/i) \; .
$$
We have  $p_i(m) \geq \frac{mA+B}{i}  - \frac{C_0}{i^2}$. It gives us
$$
C_1 (Am+B) \ge \frac{C_1 C_0}{i} + \frac{C_3}{Am+B+1} + C_2 \;.
$$
This equality holds for large $i$ and $C_1$.
This completes the proof for $d=m$.

Remind that the proof is by induction on $d$ and $i$.
Consider $d > m$ and assume that we can prove the theorem for all smaller degrees.
Now we use induction on $i$.

We have $N_i(d) \leq \frac{2mi}{d}$,
therefore $N_i(d) = O \left(i c(m,d)d^{1/A} \right)$.
In particular,
for $i < 2 \,C_7\, d^2$, where the constant $C_7$ depends only on the parameters of the model and will be defined later, we have
%$i = O\left(d^2\right)$ we have
$\E N_i(d) = c(m,d) \left(i + \theta\left(Cd^{2+1/A}\right)\right)$ with some $C$.
Assume that
$$
\E N_i(d) = c(m,d) \left(i + \theta\left(C d^{2+1/A}\right)\right) \; .
$$
From (\ref{d>m}) we obtain
\begin{multline*}
\E N_{i+1}(d) =
\E N_i(d)\left(1 - p_i(d)\right) + \E N_i(d-1) p_i^1(d-1) +
 \\ +\sum_{j=2}^m \E N_i(d-j) p_i^j(d-j) + O(p_i) =
\\
= c(m,d) \left(i + \theta\left(C d^{2+1/A}\right)\right) \left(1 - p_i(d)\right) + \\
+ c(m,d-1) \left(i + \theta\left(C (d-1)^{2+1/A}\right)\right) p_i^1(d-1)
+ \theta\left(\frac{C_4 c(m,d) d^2 i d^{1/A}}{i^2}\right) =
\\
= c(m,d)(i+1) + c(m,d-1) i p_i^1(d-1) - \\ - c(m,d) i p_i(d) - c(m,d)
+ c(m,d)\theta\left(C d^{2+1/A}\right)\left(1 - p_i(d)\right) + \\
+ \frac{c(m,d)(Ad+B+1)}{Ad-A+B}\theta\left(C (d-1)^{2+1/A}\right)p_i^1(d-1) +  \theta\left(\frac{C_4 c(m,d) d^2  d^{1/A}}{i}\right)  =
\\
= c(m,d)(i+1) +  c(m,d)\theta\left(C d^{2+1/A}\right)\left(1 - p_i(d)\right) +
\\
+\frac{c(m,d)(Ad+B+1)}{Ad-A+B}\theta\left(C (d-1)^{2+1/A}\right)p_i^1(d-1) +
\theta\left(\frac{C_5 c(m,d) d^2  d^{1/A}}{i}\right) \;.
%= c(m,d)(i+1) + c(m,d)\theta\left(C d^\delta\right)\left(1 - p_i(d)\right) %+ \theta\left(\frac{ C_5 c(m,d) d^2}{i}\right)
%+ \\
%+ \frac{c(m,d)(Ad+B+1)}{Ad-A+B}\theta\left(C (d-1)^\delta\right)p_i^1(d-1) +
% \theta\left(\frac{C_4 c(m,d)}{i}\right) \;.
\end{multline*}
We need to prove that there exists a constant $C$ that
$$
C d^{2+1/A} p_i(d) \ge \frac{C(Ad+B+1)}{Ad-A+B}  (d-1)^{2+1/A} p_i^1(d-1) + \frac{C_5 d^{2+1/A}}{i}\;,
$$
\begin{multline*}
C d^{2+1/A} p_i(d) \ge  \frac{C(Ad+B+1)}{Ad-A+B}
\left(d^{2+1/A}-(2+1/A) d^{1+1/A} + C_6 d^{1/A}\right)\cdot \\ \cdot p_i^1(d-1) + \frac{C_5 d^{2+1/A}}{i}\;,
\end{multline*}
\begin{multline*}
\frac{Cd^{2+1/A}}{i}
%(1+\frac{A+1}{Ad-A+B})
\left(2A + \frac{(B-A)(2A+1)}{Ad} + O\left(\frac{d}{i^2}\right)\right)
 \ge C d^{2+1/A}  O\left(\frac{d^2}{i^2}\right)
+ \\ + \frac{C(Ad+B+1)}{Ad-A+B} C_6 d^{1/A}  \left( A \frac{d-1}{i} + B\frac{1}{i} + O\left(\frac{d^2}{i^2}\right)\right) + \frac{C_5 d^{2+1/A}}{i}\;,
\end{multline*}
$$
\frac{C d^{2+1/A}}{i} \ge
\frac{ C_7 C d^{4+1/A}}{i^2}
+ \frac{ C_8 C d^{1+1/A}}{i}
+ \frac{C_9 d^{2+1/A}}{i}\;.
$$
%$$
%C d^\delta p_i(d) \ge \frac{C(Ad+B+1)p_i(d-1)}{Ad-A+B} \left(d^{1/A} - \frac{1}{A}d^{1/A-1} +
%\frac{1-A}{A^2}d^{1/A-2} + \theta\left(C_5 d^{1/A-3}\right) \right) + \frac{C_4}{i}\;,
%$$
%$$
%\frac{C(B+A)}{A}d^{1/A-1} \ge  C C_6 d^{1/A-2} + C_4  \;.%d^{\frac{1}{A}} + \frac{C_5 C d^\frac{1}{A}  d^2}{i},
%$$
This inequality holds for large $C \geq C_{10}$ and $d \geq d_1$.
%\textbf{[Seems that some explanation is needed to explain that we can beat the first term. Actually, larger $C$ allows to cover more $d$ at the beginning, so we can make the fraction $d/i^2$ as small as we need.]}
For constant number of small $d<d_1$ there exists a function $f(d)>0$ such that
$$
f(d) d^{2+1/A}p_i(d) \ge f(d-1)\frac{Ad+B+1}{Ad-A+B}(d-1)^{2+1/A}p_i^1(d-1) + \frac{C_5 d^{2+1/A}}{i}\;.
$$
Thus the final $C$ is $\max \left\{ C_{10} , \max_{d<d_1} \{f(d)\}  \right\}$.
This concludes the proof.

\subsection*{Proof of Theorem \ref{Concentration}}

To prove Theorem \ref{Concentration} we need the Azuma--Hoeffding inequality:

\begin{theorem}[Azuma, Hoeffding]\label{Azuma}
Let $(X_i)_{i=0}^{n}$ be a martingale such that
$|X_i - X_{i-1}| \le c_i$ for any $1 \le i \le n$.
Then
$$
\Prob\left(|X_n - X_0| \ge x \right) \le 2e^{-\frac{x^2}{2\sum_{i=1}^n c_i^2}}
$$
for any $x>0$.
\end{theorem}

Suppose we are given some $\delta > 0$.
Fix $n$ and $d$: $1 \le d \le n^{\frac{A-\delta}{4A+2}}$.
Consider the random variables
$X_i(d) = \E(N_n(d)\mid G_m^i)$,
$i = 0, \ldots, n$.

%Let us explain the notation $\E(N_n(d)\mid G_m^i)$.
%Denote by $\mathfrak{G}_m^n$ the probability space of graphs we obtain after $n$-th step of the process.
%We construct the graph $G_m^n \in \mathfrak{G}_m^n$ by induction.
%For any $t \le n$ there exists a unique
%$G_m^t \in \mathfrak{G}_m^t$ such that
%$G_m^n$ is obtained from $G_m^t$.
%So $\E(N_n(d)\mid G_m^t)$ is the expectation of the number of vertices with degree $d$ in $G_m^n$
%if at the step $t$ we have the graph $G_m^t$.
%Note that
%$X_0(d) = \E N_n(d)$ and $X_n(d) = N_n(d)$.
%From the definition of $G_m^n$ it follows that $X_n(d)$ is a martingale.

Let us explain the meaning of the random variable $\E(N_n(d)\mid G_m^i)$.
%Consider the evolution of graph $G_m^n$.
%For any $t \le n$ there exists a unique $G_m^t$ ancestor of $G_m^n$ (the induced subgraph of $G_m^n$ on the first $t$ vertices).
For any $t \le n$ let $\E(N_n(d)\mid G_m^t)$ be the expectation of the number of vertices with degree $d$ we may have at the step $n$ of the process $G_m^t$
if we fix first $t$ steps of the evolution and allow the rest $n-t$ steps to be arbitrary.
Note that $X_0(d) = \E N_n(d)$ and $X_n(d) = N_n(d)$.
It is easy to see that $X_n(d)$ is a martingale.

We will prove below that for any $i = 0, \ldots, n-1$
$$
|X_{i+1}(d) - X_{i}(d)| \le M d,
$$
where $M > 0$ is some constant.
Theorem follows from this statement immediately.
Put $c_i = M d$ for all $i$.
Then from Azuma--Hoeffding inequality it follows that
$$
\Prob\left(|N_n(d) - \E N_n(d)| \ge d \, \sqrt{n}  \, \log{n}\right) \le
2\exp\left\{-\frac{n \,d^2 \,\log^2{n} }{2\,n\,M^2 d^2 } \right\} = O\left(n^{-\log{n}}\right).
$$
If $d \le n^{\frac{A-\delta}{4A+2}}$, then the value of
$\frac{n}{d^{1+1/A}}$ is considerably greater than $d \, \log{n}\, \sqrt{n}$.
%This means that we have
%$\left(k \, \sqrt{n}\, \log^2{n}\right)/\left(n/k^2  \right) = o(1)$.
This is exactly what we need.

It remains to estimate the quantity $|X_{i+1}(d) - X_{i}(d)|$.
The proof is by a direct calculation.

Fix $0 \le i \le n-1$ and some graph $G_m^{i}$.
Note that
\begin{multline*}
\left|\E\left(N_n(d)\mid G_m^{i+1}\right) - \E\left(N_n(d)\mid G_m^{i}\right)\right| \le \\
\le \max_{\tilde G_m^{i+1}\supset G_m^{i}}  \left\{\E\left(N_n(d)\mid \tilde G_m^{i+1}\right)\right\} -
\min_{\tilde G_m^{i+1}\supset G_m^{i}}  \left\{\E\left(N_n(d)\mid \tilde G_m^{i+1}\right)\right\}.
\end{multline*}

Put
$\hat G_m^{i+1} = \arg \max \E(N_n(d)\mid \tilde G_m^{i+1})$,
$\bar G_m^{i+1} = \arg \min \E(N_n(d)\mid \tilde G_m^{i+1})$.
We need to estimate the difference
$\E(N_n(d)\mid \hat G_m^{i+1}) - \E(N_n(d)\mid \bar G_m^{i+1})$.

For $i+1 \le t \le n$ put
$$
\delta_t^i(d) = \E(N_t(d)\mid \hat G_m^{i+1}) - \E(N_t(d)\mid \bar G_m^{i+1}).
$$

First let us note that for $t \le C_{11} d^2$,
%where the constant $C_{11}$ will be defined later,
then we have $\delta_t^i(d) \le \frac{2mt}{d} \le Md$ for some constant $M$.

Now we want to prove that $\delta^i_n(d) \le M d$ by induction.
Suppose that $n=i+1$.
Fix $G_m^{i}$.
Graphs $\hat G_m^{i+1}$ and $\bar G_m^{i+1}$ are obtained from the graph $G_m^{i}$ by
adding the vertex $i+1$ and $m$ edges. Therefore
$\delta^i_{i+1}(d) \le 2m.$

Now consider $t$: $i \le t \le n-1$, $t > C_{11} d^2$.
Note that
$$
\E \left(N_{t+1} (m) \mid G_m^i\right) =
\E\left(N_{t} (m) \mid G_m^i\right)\left(1 - p_t(m)\right) + 1 + O(1/t) \;,
$$
$$
\E\left(N_{t+1} (d) \mid G_m^i\right) =
\E\left(N_{t} (d) \mid G_m^i\right)\left(1 - p_t(d)\right) +
$$
$$
+ \E\left(N_{t} (d-1) \mid G_m^i\right)p_t^1(d-1)
+ \sum_{j=2}^m \E\left(N_{t} (d-j) \mid G_m^i\right)p_t^j(d-j) + O(1/t)
, \,\,\, d \ge m+1 \;.
$$
%$$
%\E\left(N_{t+1} (d) \mid G_m^i\right) =
%\E\left(N_{t} (d) \mid G_m^i\right)\left(1 - p_t(d)\right) +
%$$
%$$
%+ \E\left(N_{t} (d-1) \mid G_m^i\right)p_t^1(d-1)
% + O(1/t)
%, \,\,\, d > 2m \;.
%$$
We obtained the same equalities in the proof of Theorem \ref{Expectation}, see (\ref{d=m})-(\ref{d>m}).
Replace $G_m^i$ by $\hat G_m^i$ or $\bar G_m^i$ in these equalities.
Substracting the equalities with $\bar G_m^i$  from the equalities with $\hat G_m^i$
%and using the inequality $(a+b)I(a+b>0) \le aI(a>0) + bI(b>0)$,
we get (for $d > m$)
\begin{multline}\label{eq:delta}
\delta^i_{t+1}(d) = \delta^i_t(d) \left(1 - p_t(d) \right) + \delta^i_t(d-1)p_t^1(d-1) + O\left(\frac{\E N_t(d) d^2}{t^2}\right) + O\left( \frac{1}{t}\right)
= \\ =
\delta^i_t(d) \left(1 - p_t(d) \right) + \delta^i_t(d-1)p_t^1(d-1) + \theta\left(\frac{C_{12} d}{t}\right)
\; .
\end{multline}
Here we used that $\E N_t(d) = O\left(td^{-1-1/A} + d\right) = O(t/d)$.
From this recurrent relation it is easy to obtain by induction that
$\delta^i_n(d) \le M d$ for some $M$.
\begin{multline*}
\delta^i_{t+1}(d) \le
Md \left(1 - p_t(d) \right) + M(d-1)p_t^1(d-1) + \frac{C_{12} d}{t} \le \\
\le Md - \frac{MA(2d-1)}{t} -\frac{MB}{d} + \frac{C_{13}Md^3}{t^2} + \frac{C_{12} d}{t} \le Md
\;
\end{multline*}
for sufficiently large $M$.
This concludes the proof of Theorem \ref{Concentration}.

\subsection*{Proof of Theorem \ref{P_2}}

Let us give the sketch of the proof of Theorem \ref{P_2}.
We can prove this theorem by induction.
Note that
$$
P_2(n) = \sum_{d = m}^\infty N_n(d)\frac{d(d-1)}{2} \; .
$$
Therefore
\begin{multline*}
\E P_2(i+1) = \sum_{d = m}^\infty \E N_{i+1}(d)\frac{d(d-1)}{2} =
\E P_2(i) + \frac{m(m-1)}{2} + \sum_{d = m}^\infty \E N_i(d) p_i(d) d \sim \\
\sim
\E P_2(i) + \frac{m(m-1)}{2} + \sum_{d = m}^\infty \frac{(Ad+B) d \E N_i(d)}{i}
=\E P_2(i) \left(1 + \frac{2A}{i}\right) + \frac{m(m-1)}{2} + \\ + \sum_{d = m}^\infty \frac{(A+B) d \E N_i(d)}{i}
= \E P_2(i) \left(1 + \frac{2A}{i}\right) + 2m(A+B) + \frac{m(m-1)}{2}.
\end{multline*}
So we obtain
\begin{multline*}
\E P_2(n) \sim \left(2m(A+B) + \frac{m(m-1)}{2}\right) \sum_{t=1}^n \prod_{i=t+1}^n \left(1 + \frac{2A}{i}\right) \sim \\
\sim \left(2m(A+B) + \frac{m(m-1)}{2}\right) \sum_{t=1}^n \frac{n^{2A}}{t^{2A}} \; .
\end{multline*}

If $2A<1$ then
$$
\E P_2(n) \sim \left(2m(A+B) + \frac{m(m-1)}{2}\right) \frac{n}{1-2A} \;.
$$
If $2A=1$ then
$$
\E P_2(n) \sim \left(2m(A+B) + \frac{m(m-1)}{2}\right) n \log(n) \;.
$$
If $2A>1$ then
$$
\E P_2(n) = O\left( n^{2A}\right)\;.
$$

Note that if $2A \le 1$, then the structure of an arbitrary graph $G_m^{n_0}$ does not affect the asymptotic of $\E P_2(n)$. If $2A>1$, then $G_m^{n_0}$ affects only the constant in $O\left( n^{2A}\right)$.

We computed the expectation of $P_2$.
One can prove concentration using standard martingale methods, although the proof is not trivial in this case. Here we need the fact that the maximum degree $\Delta_n$ grows as $n^{A}$, which can be shown using an induction. Let us consider the case $1-2A>0$.
The intuition behind this proof is the following. If we draw an edge to some vertex then this edge increase the expected final degree of this vertex by $(n/i)^A$. Finally, the expected number of $P_2$ increases by at most $n^{A}(n/i)^A = n^{2A}/i^A$ (we multiply the number of extra edges by the maximum possible degree of a vertex). Now, the sum of the squares of these values (see $\sum_{i=1}^n c_i^2$ in Theorem~\ref{Azuma}) is of order $n^{1+2A}$.
So, in Azuma's inequality we can take $x$ growing faster than $n^{1/2+A}$. Note that in this case $x$ can be taken smaller than $\E P_2(n)$ which gives concentration. In the case $1-2A<0$ we are not able to get concentration, but it is possible to get asymptotic from Theorem~\ref{P_2}.


\begin{thebibliography}{24}

\bibitem{BA_Review} R. Albert, A.-L. Barab\'asi, Statistical mechanics of complex networks, Reviews of modern physics, vol. 74, pp. 47--97 (2002)

\bibitem{BioInfoPrior} S. Bansal, S. Khandelwal, L.A. Meyers, Exploring biological network structure with clustered random networks, BMC Bioinformatics, 10:405~(2009)

\bibitem{BA} A.-L. Barab\'asi, R. Albert, Science 286, 509 (1999);
A.-L. Barab\'asi, R. Albert, H. Jeong, Physica A 272, 173 (1999);
R. Albert, H. Jeong, A.-L. Barab\'asi, Nature 401, 130 (1999)

\bibitem{Generating} V. Batagelj, U. Brandes, Efficient generation of large random networks, Phys. Rev. E, vol. 71, 036113 (2005)

\bibitem{Networks} S. Boccaletti, V. Latora, Y. Moreno, M. Chavez, D.-U. Hwang, Complex networks: Structure and dynamics, Physics reports, vol. 424(45), pp. 175-308 (2006)

\bibitem{Math_Results} B. Bollob\'as, O.M. Riordan, Mathematical results on scale-free random graphs, Handbook of Graphs and Networks: From the Genome to the Internet, pp. 1-3~(2003)

\bibitem{LCD_degrees} B. Bollob\'as, O.M. Riordan, J. Spencer, G. Tusn\'ady, The degree sequence of a scale-free random graph process, Random Structures and Algorithms, vol. 18(3), pp. 279-290~(2001)

\bibitem{Chayes} C. Borgs, M. Brautbar, J. Chayes, S. Khanna, B. Lucier, The power of local information in social networks, preprint~(2012)

\bibitem{Broder} A. Broder, R. Kumar, F. Maghoul, P. Raghavan, S. Rajagopalan, R. Stata, A. Tomkins, J. Wiener, Graph structure in the web, Computer Networks, vol. 33(16), pp. 309-320~(2000)

\bibitem{Buckley_Osthus} P.G. Buckley, D. Osthus, Popularity based random graph models leading to a scale-free degree sequence, Discrete Mathematics, vol. 282, pp. 53-63~(2004)

\bibitem{CooperA1A2} C. Cooper, Distribution of Vertex Degree in Web-Graphs, Combinatorics, Probability and Computing, vol. 15 , pp 637-661,~(2006)

\bibitem{CooperFrieze} C. Cooper, A. Frieze, A General Model of Web Graphs, Random Structures and Algorithms, 22(3), pp. 311-335~(2003)

\bibitem{Pralat} C. Cooper, P. Pra\l at, Scale-free graphs of increasing degree, Random Structures and Algorithms, vol. 38(4), pp. 396-421,~(2011)

\bibitem{Ark} M. Deijfen, H. van den Esker, R. van der Hofstad, G. Hooghiemstra, A preferential attachment model with random initial degrees, Ark. Mat., vol.47, pp. 41-72~(2009)

\bibitem{Mori_Clustering} N. Eggemann, S.D. Noble, The clustering coefficient of a scale-free random graph, Discrete Applied Mathematics, vol. 159(10), pp. 953-965~(2011)

\bibitem{F-F-F} M. Faloutsos, P. Faloutsos, Ch. Faloutsos, On power-law relationships of the Internet topology, Proc. SIGCOMM'99~(1999)

\bibitem{Gr} E.A. Grechnikov, An estimate for the number of edges between vertices of given degrees in random graphs in the Bollob\'as--Riordan model, Moscow Journal of Combinatorics and Number Theory, vol.1(2), pp. 40--73~(2011)

\bibitem{Holme_Kim} P. Holme, B.J. Kim, Growing scale-free networks with tunable clustering, Phys. Rev. E, vol. 65(2), 026107~(2002)

\bibitem{Mori} T.F. M\'ori, The maximum degree of the Barab\'asi-Albert random tree, Combinatorics, Probability and Computing,  vol. 14, pp. 339-348,~(2005)

\bibitem{Prior1} M.\'A. Serrano, M. Bogu\~n\'a, Tuning clustering in random networks with arbitrary degree distributions, Phys. Rev. E, vol. 72(3),036133~(2005)

\bibitem{Prior2} E. Volz, Random Networks with Tunable Degree Distribution and Clustering, Phys. Rev. E, vol. 70(5), 056115~(2004)

\bibitem{RAN} T. Zhou, G. Yan and B.-H. Wang, Maximal planar networks with large clustering coefficient and power-law degree distribution journal, Phys. Rev. E, vol. 71(4), 46141~(2005)

\newpage

\end{thebibliography}
\end{document}